\documentclass[12pt]{amsart}

\newtheorem{THM}{Theorem}[section]
\newtheorem{LMA}[THM]{Lemma}
\newtheorem{PROP}[THM]{Proposition}

\numberwithin{equation}{section}

\usepackage{amsmath}
\usepackage{amssymb}
\usepackage{euscript}

\newcommand{\showon}{\begin{eqnarray*}}
\newcommand{\showoff}{\end{eqnarray*}}

\newcommand{\none}{\varnothing}
\newcommand{\drop}{\smallsetminus}

\newcommand{\one}{\boldsymbol{1}}
\newcommand{\zero}{\boldsymbol{0}}

\newcommand{\A}{\EuScript{A}}

\newcommand{\J}{\EuScript{J}}
\newcommand{\K}{\EuScript{K}}
\renewcommand{\L}{\EuScript{L}}
\newcommand{\M}{\EuScript{M}} 
\newcommand{\N}{\EuScript{N}}

 \newcommand{\RR}{\mathbb{R}}
 \renewcommand{\S}{\EuScript{S}}

\newcommand{\w}{\mathbf{w}}
 \newcommand{\x}{\mathbf{x}}

\newcommand{\Y}{\EuScript{Y}} \newcommand{\y}{\mathbf{y}}

\newcommand{\z}{\mathbf{z}}

\begin{document}

\title{Rank three matroids are Rayleigh}

\author{David G. Wagner}
\address{Department of Combinatorics and Optimization\\
University of Waterloo\\
Waterloo, Ontario, Canada\ \ N2L 3G1}
\email{\texttt{dgwagner@math.uwaterloo.ca}}


\begin{abstract}
A Rayleigh matroid is one which satisfies a set of inequalities
analogous to the Rayleigh monotonicity property of linear
resistive electrical networks.  We show that every matroid of
rank three satisfies these inequalities.
\end{abstract}

\maketitle

\section{Introduction.}

(For the basic concepts of matroid theory we refer the reader
to Oxley's book \cite{Ox}.)

A linear resistive electrical network can be represented as a graph
$G=(V,E)$ together with a set of positive real numbers $\y=\{y_{e}:\
e\in E\}$ that specify the conductances of the
corresponding elements.  In 1847 Kirchhoff \cite{Ki}
determined the effective
conductance of the network measured between vertices $a,b\in V$
as a rational function $\Y_{ab}(G;\y)$ of the conductances $\y$.
This formula can be generalized directly to any matroid.

For electrical networks the following property is physically 
intuitive:\ if $y_{c}>0$ for all $c\in E$ then for any $e\in E$,
$$\frac{\partial}{\partial y_{e}}\Y_{ab}(G;\y)\geq 0.$$
That is, by increasing the conductance of the element $e$ we can not
decrease the effective conductance of the network as a whole.
This is known as the \emph{Rayleigh monotonicity} property.

Informally, a matroid has the Rayleigh property if it satisfies
inequalities analogous to the Rayleigh monotonicity property of
linear resistive electrical networks. 
While there are non--Rayleigh matroids of rank four or more, we
show here that every matroid of rank (at most) three is Rayleigh,
answering a question left open by Choe and Wagner \cite{CW}.

Let $\M$ be a matroid with ground--set $E$, and fix indeterminates
$\y:=\{y_{e}:\ e\in E\}$ indexed by $E$.  For a basis $B$ of $\M$ let
$\y^{B}:=\prod_{e\in B}y_{e}$, and let $M(\y):=\sum_{B\in\M}\y^{B}$ with
the sum over all bases of $\M$.  Since $M(\y)$ is insensitive to the
presence of loops we generally consider only loopless matroids, and
regard $\M$ as its set of bases.

For disjoint subsets $I,J$ of $E$, let $\M_{I}^{J}$ denote the minor
of $\M$  obtained by contracting $I$ and deleting $J$.  We use the 
nonstandard convention that if $I$ is dependent then $\M_{I}^{J}$ is
empty, so that in general
$$\M_{I}^{J}:=
\{B\drop I:\ B\in\M\ \mathrm{and}\ I\subseteq B\subseteq E(\M)\drop 
J\}.$$

The matroid $\M$ is a \emph{Rayleigh} matroid provided that whenever
$y_{c}>0$ for all $c\in E$, then for every pair of distinct $e,f\in E$,
$$\Delta 
M\{e,f\}(\y):=M_{e}^{f}(\y)M_{f}^{e}(\y)-M_{ef}(\y)M^{ef}(\y)\geq 0.$$
See Section 3 of Choe and Wagner \cite{CW} for more detailed motivation 
of this definition.  Rayleigh matroids are ``balanced'' in the sense of
Feder and Mihail \cite{FM}, and for binary matroids these conditions
are equivalent.  For example, every sixth--root of unity matroid -- 
in particular every regular matroid -- is Rayleigh (Proposition 5.1 and
Corollary 4.9 of \cite{CW}).  Since graphic matroids are regular this
generalizes the physical assertion that linear resistive electrical networks
satisfy Rayleigh monotonicity.  One of the main questions left open in
\cite{CW} is whether or not every matroid of rank three is Rayleigh.
Here we show that this is indeed the case.

\begin{THM} Every matroid of rank three is Rayleigh.
\end{THM}
In contrast to this theorem there are several matroids of rank
four which are known not to be Rayleigh, among them the matroids
$\S_{8}$ and $\J'$ discussed in \cite{CW}.

As a concrete but fairly representative consequence of Theorem 1.1,
let $E$ be a finite non--collinear set of points in a projective
plane, and let $\M$ be the set of non--collinear unordered triples
of points in $E$.  Assign a positive real number $y_{c}$ to each
$c\in E$, and consider the probability space $\Omega(\M,\y)$ which
assigns to each $B\in\M$ the probability $\y^{B}/M(\y)$.  Since
$\M$ is a rank three matroid it is Rayleigh, by Theorem 1.1.  A
short calculation shows that for distinct $e,f\in E$:
$$\frac{M_{ef}(\y)}{M_{e}(\y)}\leq\frac{M_{f}(\y)}{M(\y)}.$$
That is, in $\Omega(\M,\y)$ the probability that a random basis
$B\in\M$ contains $f$, given that it contains $e$, is at most
the probability that a random basis contains $f$.  In short, the
events $e\in B$ and $f\in B$ are \emph{negatively correlated}
for any distict $e,f\in E$.  This probabilistic point of view
is carried further by Feder and Mihail \cite{FM} and Lyons
\cite{Ly}.

Several conversations and correspondences with Jim Geelen, Sandra Kingan,
and Bruce Reznick helped to clarify my thoughts on this problem, for 
which I thank them sincerely.

\section{Preliminaries.}

To simplify notation, when calculating with Rayleigh matroids we
will henceforth usually omit reference to the variables $\y$ -- writing
$M_{I}^{J}$ instead of $M_{I}^{J}(\y)$ \emph{et cetera} -- unless
a particular substitution of variables requires emphasis.  We will
also write ``$\y>\zero$'' as shorthand for ``$y_{c}>0$ for all $c\in 
E$''.

We require the following facts from \cite{CW}.
\begin{PROP}[Section 3 of \cite{CW}]
The class of Rayleigh matroids is closed by taking duals and minors.
\end{PROP}
\begin{proof}[Sketch of proof]
For the matroid $\M^{*}$ dual to $\M$ and for $e,f\in E(\M^{*})$,
$$\Delta M^{*}\{e,f\}(\y)=\y^{2E}\Delta M\{e,f\}(\one/\y)$$
in which $\one/\y:=\{1/y_{c}:\ c\in E\}$.
From this it follows that $\M^{*}$ is Rayleigh if $\M$ is.

For distinct $e,f,g\in E(\M)$,
$$\Delta M^{g}\{e,f\}=\lim_{y_{g}\rightarrow 0}\Delta M\{e,f\}$$
and
$$\Delta M_{g}\{e,f\}=\lim_{y_{g}\rightarrow\infty}\frac{1}{y_{g}^{2}}
\Delta M\{e,f\}.$$
From this it follows that if $\M$ is Rayleigh then the deletion
$\M^{g}$ and the contraction $\M_{g}$ are also Rayleigh.  The case of a
general minor follows by iteration of these two cases.
\end{proof}
(The class of Rayleigh matroids is also closed by $2$-sums, but we 
will not use this fact.)

For polynomials $A(\y)$ and $B(\y)$ in $\RR[y_{c}:\ c\in E]$,
we write $A(\y) \gg  B(\y)$ to mean that every coefficient of
$A(\y)-B(\y)$ is nonnegative.  Certainly, if $A(\y)\gg 0$ then
$A(\y)\geq 0$ for all $\y>\zero$, but not conversely.
Making the substitution $y_{c}= x_{c}^{2}$ for each $c\in E$,
we have $A(\y)\geq 0$ for all $\y>\zero$ if and only if
$A(\x^{2})\geq 0$ for all $\x \in\RR^{E}$; such a form $A(\x^{2})$
is said to be \emph{positive semidefinite}.
Artin's solution to Hilbert's 17th problem asserts that every
positive semidefinite form can be written as a positive sum of
squares of rational functions, but the proof is nonconstructive.
Reznick \cite{R} gives an excellent survey of Hilbert's 17th
problem.  To prove Theorem 1.1 we will write $\Delta M\{e,f\}(\y)$
as a positive sum of monomials and squares of polynomials in $\y$.

Regarding the Rayleigh property, one may restrict attention to the
class of simple matroids (although it is not always useful to do so)
for the following reason.
We may assume that $\M$ is loopless, as remarked above.
If $a, a_{1},\ldots,a_{k}$ are parallel elements in $\M$, then let
$\N$ be obtained from $\M$ by deleting $a_{1},\ldots,a_{k}$.
Letting $w_{c}:=y_{c}$ if $c\in E(\N)\drop\{a\}$ and
$w_{a}:=y_{a}+y_{a_{1}}+\cdots+y_{a_{k}}$, one sees that 
$M(\y)=N(\w)$.   A little calculation shows that $\M$ is Rayleigh if
and only if $\N$ is Rayleigh.  Repeating this reduction as required,
we find a simple matroid $\L$ and a substitution of variables $\z=\z(\y)$
such that $M(\y)=L(\z)$, and such that $\M$ is Rayleigh if and only if
$\L$ is Rayleigh.

It is very easy to see that matroids of rank one or two are Rayleigh.
\begin{PROP}
If $\M$ has rank at most two then $\Delta M\{e,f\}\gg 0$ for all
distinct $e,f\in E(\M)$.  Consequently, $\M$ is Rayleigh.
\end{PROP}
\begin{proof}
By the above remarks, we may assume that $\M$ is simple.
Let the ground--set of $\M$ be $E=\{1,2,\ldots,m\}$.

If $\M$ has rank one then $M(\y)=y_{1}+y_{2}+\cdots+y_{m}$, so
$M_{ef}=0$ for all distinct $e,f\in E$, and hence
$\Delta M\{e,f\}=M_{e}^{f}M_{f}^{e}\gg 0$.

If $\M$ has rank two then $M(\y)=\sum_{1\leq i<j\leq m} y_{i}y_{j}$ is
the second elementary symmetric function of $\y$.  By symmetry
we only need to show that $\Delta M\{1,2\}\gg 0$.  
Since $M_{1}^{2}=M_{2}^{1}=y_{3}+y_{4}+\cdots +y_{m}$ and
$M_{12}=1$ and
$$M^{12}=\sum_{3\leq i<j\leq m}y_{i}y_{j},$$
it follows that
$$\Delta M\{1,2\}=\sum_{3\leq i\leq j\leq m}y_{i}y_{i},$$
proving that $\Delta M\{1,2\}\gg 0$.
\end{proof}

The case of rank three matroids is much more interesting -- the 
polynomial $\Delta M\{e,f\}$ can have terms with negative 
coefficients, as happens already for the graphic matroid $\K$
of the complete graph $\mathsf{K}_{4}$ on four vertices.  With the
ground--set of $\K$ labelled as in Figure 3(IV), we have
$$\Delta K\{1,2\}=(y_{3}y_{4}-y_{5}y_{6})^{2}.$$
As will be seen in Table 3, however, in some sense this is the worst
that can happen in rank three.

\section{A reduction lemma for any rank.}

For distinct elements $e,f,g\in E(\M)$, a short calculation shows
that
$$\Delta M\{e,f\}=
y_{g}^{2}\Delta M_{g}\{e,f\}+y_{g}\Theta M\{e,f|g\}+\Delta 
M^{g}\{e,f\}$$
in which
\showon
\Delta M_{g}\{e,f\} &=& M_{eg}^{f}M_{fg}^{e}-M_{efg}M_{g}^{ef},\\
\Delta M^{g}\{e,f\} &=& M_{e}^{fg}M_{f}^{eg}-M_{ef}^{g}M^{efg},
\showoff
and
the \emph{central term for $\{e,f\}$ and $g$ in $\M$} is defined by
$$\Theta M\{e,f|g\}:=M_{e}^{fg}M_{fg}^{e}+M_{f}^{eg}M_{eg}^{f}
-M_{g}^{ef}M_{ef}^{g}-M_{efg}M^{efg}.$$

For a subset $S$ of $E(\M)$, we use $\overline{S}$ to denote the
closure of $S$ in $\M$.

\begin{LMA}
Let $\M$ be a matroid, and let $e,f,g\in E(\M)$ be
distinct elements. If $\{e,f,g\}$ is dependent in $\M$
then $\Theta M\{e,f|g\}\gg 0$.
\end{LMA}
\begin{proof}
To prove this we exhibit an injective function
$$\left(\M_{g}^{ef}\times \M_{ef}^{g}\right) \cup
\left(\M_{efg}\times \M^{efg}\right)
\longrightarrow 
\left(\M_{e}^{fg} \times \M_{fg}^{e}\right) \cup
\left(\M_{f}^{eg}\times \M_{eg}^{f}\right)$$
such that if $(B_{1},B_{2})\mapsto(A_{1},A_{2})$ then
$\y^{A_{1}}\y^{A_{2}}=\y^{B_{1}}\y^{B_{2}}$.

Since $\{e,f,g\}$ is dependent it follows that $\M_{efg}=\none$,
so let $B_{1}\in\M_{g}^{ef}$ and $B_{2}\in\M_{ef}^{g}$.  Let
$L:=\overline{B_{1}\drop\{g\}}$.  We claim that either
$e\not\in L$ or $f\not\in L$.  To see this, suppose not -- then 
$g\in\overline{\{e,f\}}\subseteq L$, which contradicts the fact that
$B_{1}$ is a basis.  If $e\not\in L$ then let
$A_{1}:=B_{1}\cup\{e\}\drop\{g\}$ and $A_{2}:=B_{2}\cup\{g\}\drop\{e\}$.
If $e\in L$ then $f\not\in L$, so let
$A_{1}:=B_{1}\cup\{f\}\drop\{g\}$ and $A_{2}:=B_{2}\cup\{g\}\drop\{f\}$.
It is easy to see that in either case both
$A_{1}$ and $A_{2}$ are bases of $\M$.

Notice that for $(A_{1},A_{2})$ in the image of this function, 
$A_{1}\in\M_{e}^{fg}\cup\M_{f}^{eg}$ and this union is disjoint.
If $\A_{1}\in\M_{e}^{fg}$ then let $B'_{1}:=A_{1}\cup\{g\}\drop\{e\}$
and $B'_{2}:=A_{2}\cup\{e\}\drop\{g\}$, while if $\A_{1}\in\M_{f}^{eg}$
then let $B'_{1}:=A_{1}\cup\{g\}\drop\{f\}$
and $B'_{2}:=A_{2}\cup\{f\}\drop\{g\}$.  In either case we have
$(B'_{1},B'_{2})=(B_{1},B_{2})$ showing that the function
$(B_{1},B_{2})\mapsto(A_{1},A_{2})$ is injective.

This construction provides the desired weight--preserving injection.
\end{proof}

Lemma 3.1 has the following consequence which might be helpful in
the investigation of Rayleigh matroids of rank four or more.

\begin{PROP}
Let $\M$ be a minor--minimal non--Rayleigh matroid, and let $e,f\in 
E(\M)$ and $\y>\zero$ be such that $\Delta M\{e,f\}<0$.  Then
$\{e,f\}$ is closed in $\M$.
\end{PROP}
\begin{proof}
If $g\in E(\M)\drop\{e,f\}$ is such that $\{e,f,g\}$ is dependent, then
$\Theta M\{e,f|g\}\gg 0$ by Lemma 3.1.  From this it follows that if
$\y>\zero$ then
$$\Delta M\{e,f\}=y_{g}^{2}\Delta M_{g}\{e,f\}+
y_{g}\Theta M\{e,f|g\}+\Delta M^{g}\{e,f\}\geq 0,$$
since every proper minor of $\M$ is Rayleigh.  As this contradicts
the hypothesis we conclude that $\{e,f\}$ is closed in $\M$.
\end{proof}

The following consequence of Lemma 3.1 is relevant to the present purpose.

\begin{LMA}
Let $\M$ be a matroid of rank three, and let $e,f\in E(\M)$.
If $g\in E(\M)\drop\{e,f\}$ is such that $\{e,f,g\}$ is dependent in
$\M$ then $\Delta M\{e,f\}(\y)\gg \Delta M^{g}\{e,f\}(\y)$.
\end{LMA}
\begin{proof}
Since
$$\Delta M\{e,f\}-\Delta M^{g}\{e,f\}=
y_{g}^{2}\Delta M_{g}\{e,f\}+y_{g}\Theta M\{e,f|g\},$$
the inequality follows directly from Proposition 2.2 and Lemma 3.1.
\end{proof}

\section{Matroids of rank three.}

The proof of Theorem 1.1 is completed by means of the following
\emph{Ansatz}, which was found mainly by trial and error.

For $a\in E(\M)\drop\{e,f\}$ let $L(a,e):=\overline{\{a,e\}}\drop\{a,e\}$,
let $L(a,f):=\overline{\{a,f\}}\drop\{a,f\}$, and let $U(a):= 
E(\M)\drop(\overline{\{a,e\}}\cup\overline{\{a,f\}})$.
Define the linear polynomials $B(a):=\sum_{b\in U(a)}y_{b}$,
$C(a):=\sum_{c\in L(a,e)}y_{c}$, and $D(a):=\sum_{d\in L(a,f)}y_{d}$, and
the quartic polynomials
$$T(\M;e,f,a;\y):=\left(y_{a}B(a)-C(a)D(a)\right)^{2}$$
for each $a\in E(\M)\drop \{e,f\}$ and
$$P(\M;e,f;\y):=\frac{1}{4}\sum_{a\in E(\M)\drop\{e,f\}}T(\M;e,f,a;\y).$$

\begin{PROP}
Let $\M$ be a matroid of rank three, and let $e,f\in E(\M)$ be
distinct.  With the notation above,
$$\Delta M\{e,f\}(\y) \gg  P(\M;e,f;\y).$$
\end{PROP}
\begin{proof}
By repeated application of Lemma 3.3, if necessary, we may assume that
$\{e,f\}$ is closed in $\M$, so we reduce to this case.

Both $\Delta:=\Delta M\{e,f\}(\y)$ and $P:=P(\M;e,f;\y)$ are homogeneous
of degree four in the indeterminates $\{y_{j}:\ j\in E(\M)\drop\{e,f\}\}$,
and the only monomials which occur with nonzero coefficient in either
of these polynomials have shape $y_{g}^{2}y_{h}^{2}$,
$y_{g}^{2}y_{h}y_{i}$, or $y_{g}y_{h}y_{i}y_{j}$.  The coefficient of
such a monomial in $\Delta$ depends only on the isomorphism 
type of the restriction $\M|\{e,f,g,h\}$, $\M|\{e,f,g,h,i\}$, or
$\M|\{e,f,g,h,i,j\}$, the positions of $e$ and $f$ in this
restriction, and, in the second case, the position of $g$ relative to
$e$ and $f$ in this restriction.  (The coefficient of of such a 
monomial in $P$ can depend on more information, as we shall see.)
Since $\{e,f\}$ is closed in $\M$, $\{e,f\}$ is also closed in any
such restriction $\N$.  The proposition is now proved by an exhaustive
case analysis of these configurations in $\M$.

\setlength{\unitlength}{1mm}
\begin{figure}
\begin{picture}(130,35)
\thicklines

\put(30,25){\circle*{2}}\put(33,24){$1$}
\put(15,10){\circle*{2}}\put(14,4){$2$}
\put(30,10){\circle*{2}}\put(29,4){$3$}
\put(45,10){\circle*{2}}\put(44,4){$4$}
\put(15,10){\line(1,0){30}}
\put(5,30){I.}

\put(100,25){\circle*{2}}\put(103,24){$1$}
\put(100,10){\circle*{2}}\put(103,9){$3$}
\put(115,25){\circle*{2}}\put(118,24){$2$}
\put(115,10){\circle*{2}}\put(118,9){$4$}
\put(90,30){II.}

\end{picture}
\caption{The four--element rank three matroids.}
\end{figure}
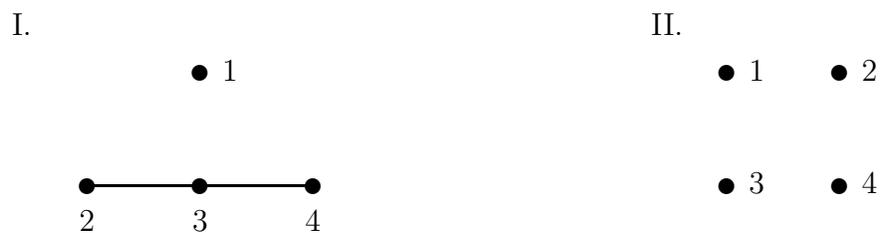

\begin{table}
$$\begin{array}{|r|l|r|c|}
\hline
\N\{e,f\} & \Delta & P & \mathrm{notes}\\ \hline
\mathrm{I}\{1,2\} & 0 - 0 = 0 & 0 &  \\
\mathrm{II}\{1,2\} & 1 - 0 = 1 & 1/2,3/4,1& \mathrm{A.}\\
\hline
\end{array}$$
\caption{Monomials of shape $y_{g}^{2}y_{h}^{2}$.}
\end{table}

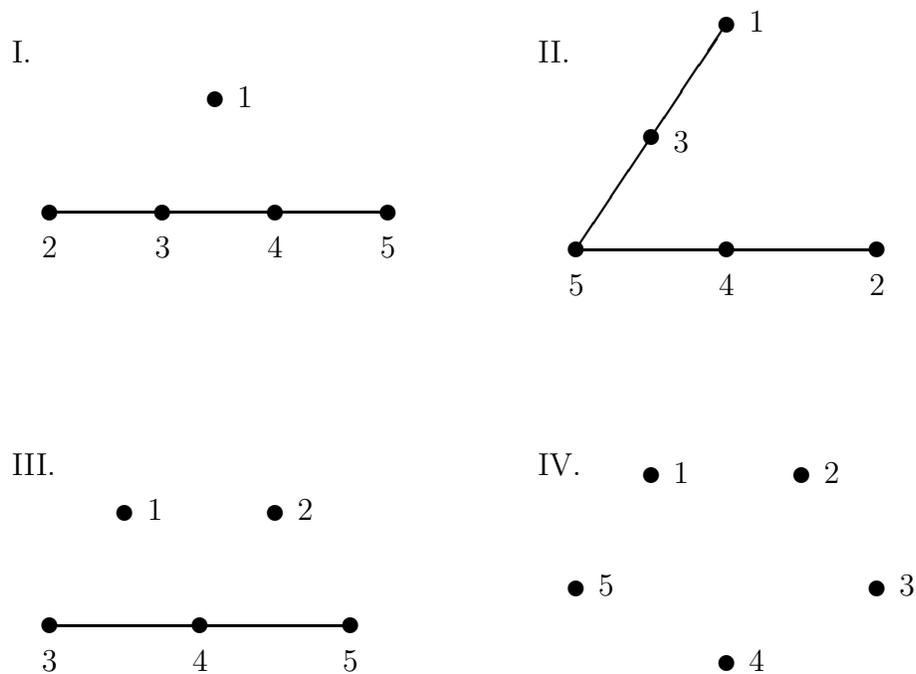
\begin{figure}
\begin{picture}(130,100)
\thicklines

\put(10,15){\line(1,0){40}}
\put(10,15){\circle*{2}}\put(9,9){$3$}
\put(30,15){\circle*{2}}\put(29,9){$4$}
\put(50,15){\circle*{2}}\put(49,9){$5$}
\put(20,30){\circle*{2}}\put(23,29){$1$}
\put(40,30){\circle*{2}}\put(43,29){$2$}
\put(5,35){III.}

\put(100,10){\circle*{2}}\put(103,9){$4$}
\put(80,20){\circle*{2}}\put(83,19){$5$}
\put(120,20){\circle*{2}}\put(123,19){$3$}
\put(90,35){\circle*{2}}\put(93,34){$1$}
\put(110,35){\circle*{2}}\put(113,34){$2$}
\put(75,35){IV.}

\put(10,70){\line(1,0){45}}
\put(10,70){\circle*{2}}\put(9,64){$2$}
\put(25,70){\circle*{2}}\put(24,64){$3$}
\put(40,70){\circle*{2}}\put(39,64){$4$}
\put(55,70){\circle*{2}}\put(54,64){$5$}
\put(32,85){\circle*{2}}\put(35,84){$1$}
\put(5,90){I.}

\put(80,65){\line(1,0){40}}
\put(80,65){\line(2,3){20}}
\put(80,65){\circle*{2}}\put(79,59){$5$}
\put(100,65){\circle*{2}}\put(99,59){$4$}
\put(120,65){\circle*{2}}\put(119,59){$2$}
\put(90,80){\circle*{2}}\put(93,78){$3$}
\put(100,95){\circle*{2}}\put(103,94){$1$}
\put(75,90){II.}

\end{picture}
\caption{The five--element rank three matroids.}
\end{figure}

\begin{table}
$$\begin{array}{|r|l|r|c|}
\hline
\N\{e,f\},g & \Delta & P & \mathrm{notes}\\ \hline
\mathrm{I}\{1,2\},3   & 0-0=0 & 0 & \\
\mathrm{II}\{1,2\},3  & 1-1=0 & 0 & \\
\mathrm{II}\{1,2\},5  & 1-1=0 & 0 & \\
\mathrm{III}\{1,2\},3 & 2-0=2 & 1/2 & \mathrm{B.}\\
\mathrm{III}\{1,3\},2 & 2-1=1 & 1/2,1 & \mathrm{C.}\\
\mathrm{III}\{1,3\},4 & 1-1=0 & 0 & \\
\mathrm{IV}\{1,2\},3  & 2-1=1 & 1/2 & \mathrm{B.}\\
\hline
\end{array}$$
\caption{Monomials of shape $y_{g}^{2}y_{h}y_{i}$.}
\end{table}

Figure 1 and Table 1 summarize the case analysis for monomials of 
shape $y_{g}^{2}y_{h}^{2}$, Figure 2 and Table 2 summarize the case
analysis for monomials of shape $y_{g}^{2}y_{h}y_{i}$, and
Figure 3 and Table 3 summarize the case analysis for monomials of
shape $y_{g}y_{h}y_{i}y_{j}$.  In each table the first column 
indicates the isomorphism class of the restriction $\N$ of $\M$,
the choice of  $\{e,f\}$ in that restriction, and (in Table 2) the
choice of $g$ in $\N$.  The second column in each table indicates the
coefficient of the relevant monomial in each term of
$$M_{e}^{f}M_{f}^{e}-M_{ef}M^{ef}=\Delta M\{e,f\},$$
respectively.  As remarked above these coefficients depend only on
$\N$, $\{e,f\}$, and $g$ and are computed from the definition by 
elementary counting.  The third column in each table indicates the
coefficient of the relevant monomial in $P:=P(\M;e,f;\y)$.  Notes in
the fourth column of each table refer to the following list of 
additional remarks regarding the coefficients of monomials of $P$:

\setlength{\unitlength}{0.8mm}
\begin{figure}
\begin{picture}(150,230)
\thicklines

\put(10,210){\line(1,0){40}}
\put(10,210){\circle*{2}}\put(9,204){$2$}
\put(20,210){\circle*{2}}\put(19,204){$3$}
\put(30,210){\circle*{2}}\put(29,204){$4$}
\put(40,210){\circle*{2}}\put(39,204){$5$}
\put(50,210){\circle*{2}}\put(49,204){$6$}
\put(30,225){\circle*{2}}\put(33,224){$1$}
\put(5,225){I.}

\put(100,205){\line(1,0){45}}
\put(100,205){\line(3,2){30}}
\put(100,205){\circle*{2}}\put(99,199){$6$}
\put(115,205){\circle*{2}}\put(114,199){$5$}
\put(130,205){\circle*{2}}\put(129,199){$4$}
\put(145,205){\circle*{2}}\put(144,199){$2$}
\put(115,215){\circle*{2}}\put(111,217){$3$}
\put(130,225){\circle*{2}}\put(133,224){$1$}
\put(95,225){II.}

\put(10,165){\line(1,0){45}}
\put(10,165){\circle*{2}}\put(9,159){$3$}
\put(25,165){\circle*{2}}\put(24,159){$4$}
\put(40,165){\circle*{2}}\put(39,159){$5$}
\put(55,165){\circle*{2}}\put(54,159){$6$}
\put(20,180){\circle*{2}}\put(23,179){$1$}
\put(45,180){\circle*{2}}\put(48,179){$2$}
\put(5,185){III.}

\put(100,155){\line(1,0){40}}
\put(100,155){\line(2,3){20}}
\put(110,170){\line(2,-1){30}}
\put(120,155){\line(0,1){30}}
\put(100,155){\circle*{2}}\put(99,149){$3$}
\put(120,155){\circle*{2}}\put(119,149){$5$}
\put(140,155){\circle*{2}}\put(139,149){$2$}
\put(120,165){\circle*{2}}\put(123,166){$4$}
\put(110,170){\circle*{2}}\put(106,172){$6$}
\put(120,185){\circle*{2}}\put(123,184){$1$}
\put(95,185){IV.}

\put(10,105){\line(1,0){40}}
\put(10,105){\line(2,3){20}}
\put(30,135){\line(2,-3){20}}
\put(10,105){\circle*{2}}\put(9,99){$2$}
\put(30,105){\circle*{2}}\put(29,99){$4$}
\put(50,105){\circle*{2}}\put(49,99){$3$}
\put(20,120){\circle*{2}}\put(16,122){$6$}
\put(40,120){\circle*{2}}\put(43,121){$5$}
\put(30,135){\circle*{2}}\put(29,138){$1$}
\put(5,135){V.}

\put(100,105){\line(1,0){40}}
\put(100,105){\line(2,3){20}}
\put(100,105){\circle*{2}}\put(99,99){$6$}
\put(120,105){\circle*{2}}\put(119,99){$5$}
\put(140,105){\circle*{2}}\put(139,99){$2$}
\put(110,120){\circle*{2}}\put(106,122){$4$}
\put(135,125){\circle*{2}}\put(138,124){$3$}
\put(120,135){\circle*{2}}\put(123,134){$1$}
\put(95,135){VI.}

\put(10,80){\line(1,0){40}}
\put(10,60){\line(1,0){40}}
\put(10,80){\circle*{2}}\put(9,74){$1$}
\put(10,60){\circle*{2}} \put(9,54){$2$}
\put(30,80){\circle*{2}}\put(29,74){$3$}
\put(30,60){\circle*{2}} \put(29,54){$5$}
\put(50,80){\circle*{2}}\put(49,74){$4$}
\put(50,60){\circle*{2}} \put(49,54){$6$}
\put(5,85){VII.}

\put(100,60){\line(1,0){40}}
\put(100,75){\circle*{2}} \put(99,69){$1$}
\put(100,60){\circle*{2}} \put(99,54){$4$}
\put(120,80){\circle*{2}}\put(119,74){$2$}
\put(120,60){\circle*{2}} \put(119,54){$5$}
\put(140,75){\circle*{2}} \put(139,69){$3$}
\put(140,60){\circle*{2}} \put(139,54){$6$}
\put(95,85){VIII.}

\put(65,40){\circle*{2}}\put(68,39){$1$}
\put(85,40){\circle*{2}}\put(88,39){$2$}
\put(95,25){\circle*{2}}\put(98,24){$3$}
\put(85,10){\circle*{2}}\put(88,9){$4$}
\put(65,10){\circle*{2}}\put(68,9){$5$}
\put(55,25){\circle*{2}}\put(58,24){$6$}
\put(50,40){IX.}

\end{picture}
\caption{The six--element rank three matroids.}
\end{figure}
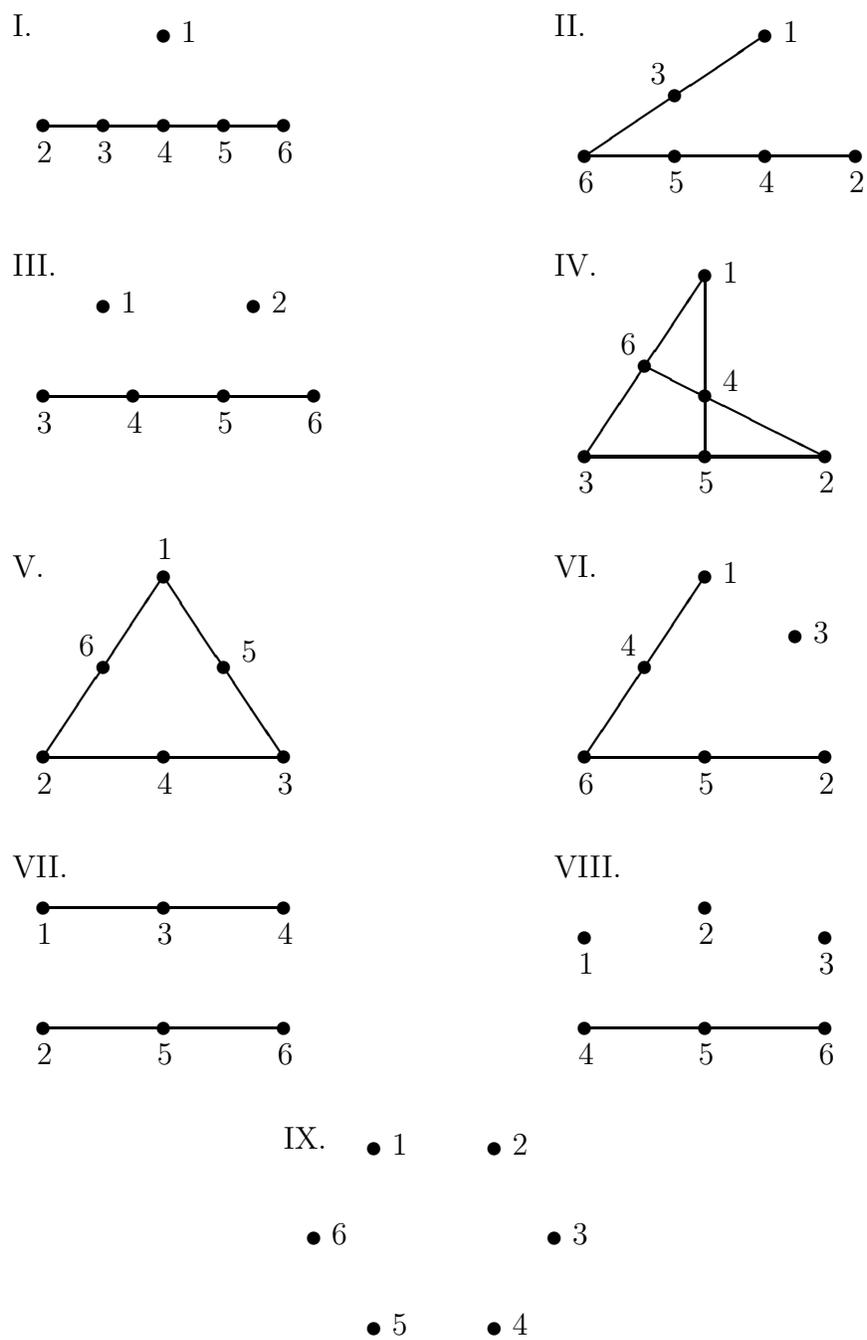

\begin{table}
$$\begin{array}{|r|l|r|c|}
\hline
\N\{e,f\} & \Delta & P & \mathrm{notes}\\ \hline
\mathrm{I}\{1,2\}    & 0-0=0  & 0 & \\
\mathrm{II}\{1,2\}   & 3-3=0  & 0 & \\
\mathrm{III}\{1,2\}  & 6-0=6  & 0 & \\
\mathrm{III}\{1,3\}  & 3-3=0  & 0 & \\
\mathrm{IV}\{1,2\}   & 2-4=-2 & -2& \mathrm{D.} \\
\mathrm{V}\{1,4\}    & 3-4=-1 & -1& \mathrm{E.} \\
\mathrm{V}\{4,5\}    & 4-3=1  & -1/2& \mathrm{F.} \\
\mathrm{VI}\{1,2\}   & 4-4=0  & -1/2& \mathrm{G.} \\
\mathrm{VI}\{1,3\}   & 5-3=2  & 0 & \\
\mathrm{VI}\{3,6\}   & 4-4=0  & 0 & \\
\mathrm{VII}\{1,2\}  & 5-4=1  & 0,1& \mathrm{H.} \\
\mathrm{VIII}\{1,2\} & 6-3=3  & 0& \\
\mathrm{VIII}\{1,4\} & 5-4=1  & 0& \\
\mathrm{IX}\{1,2\}   & 6-4=2  & 0& \\
\hline
\end{array}$$
\caption{Monomials of shape $y_{g}y_{h}y_{i}y_{j}$.}
\end{table}

\noindent
$\bullet$\ When the coefficient in the third column is zero there
is no possible location for an element $a\in E(\M)$
such that the monomial occurs in $T(\M;e,f,a;\y)$.\\
\textbf{A.}\ The monomial occurs in $T(\M;1,2,a)$ in
the term $y_{a}^{2}B(a)^{2}$ when $a=3$ or $a=4$, and in the term
$C(a)^{2}D(a)^{2}$ when $\{a\}$ is one of
$\overline{\{1,3\}}\cap\overline{\{2,4\}}$ or
$\overline{\{1,4\}}\cap\overline{\{2,3\}}$.  Either of these last two
sets might be empty instead, however.\\
\textbf{B.}\  The monomial occurs with coefficient $2$
in the term $y_{3}^{2}B(3)^{2}$ of $T(\M;1,2,3)$.\\
\textbf{C.}\ The monomial occurs with coefficient $2$
in the term $y_{2}^{2}B(2)^{2}$ of $T(\M;1,3,2)$.  If 
$\overline{\{1,2\}}\cap\overline{\{3,4\}}=\{a\}$ then
the monomial also occurs with coefficient $2$ in the term
$C(a)^{2}D(a)^{2}$ of $T(\M;1,3,a)$.  (If the above intersection is
empty then this second contribution does not occur.)\\
\textbf{D.}\ This occurs in the
term $-2y_{a}B(a)C(a)D(a)$ of $T(\M;1,2,a)$ for each 
$a\in\{3,4,5,6\}$.\\
\textbf{E.}\  This occurs in the
term $-2y_{a}B(a)C(a)D(a)$ of $T(\M;1,4,a)$ for $a=2$ and $a=3$.\\
\textbf{F.}\ This occurs in the
term $-2y_{3}B(3)C(3)D(3)$ of $T(\M;4,5,3)$.\\
\textbf{G.}\ This occurs in the
term $-2y_{6}B(6)C(6)D(6)$ of $T(\M;1,2,6)$.\\
\textbf{H.}\ If $\overline{\{1,3\}}\cap\overline{\{2,5\}}=\{a\}$ then
the monomial occurs with coefficient $4$ in the term $C(a)^{2}D(a)^{2}$
of $T(\M;1,2,a)$.  If the above intersection is empty then the monomial
does not occur.

These remarks conclude the explanation of the various coefficients of
$\Delta M\{e,f\}$ and $P(\M;e,f)$, completing the proof that
$$\Delta M\{e,f\}\gg P(\M;e,f).$$ 
\end{proof}

\begin{proof}[Proof of Theorem $1.1$]
Since $P(\M;e,f;\y)$ is a nonnegative sum of squares it follows that
$P(\M;e,f;\y)\geq 0$ for all $\y\in\RR^{E(\M)}$.  Since
$\Delta M\{e,f\}(\y) \gg  P(\M;e,f;\y)$ by Proposition 4.1
it follows that
$$\Delta M\{e,f\}(\y)\geq P(\M;e,f;\y)\geq 0$$
for all $\y>\zero$.  Hence it follows that $\M$ is Rayleigh.
\end{proof}

\end{document}